\documentclass{amsart}

\newtheorem{theorem}{Theorem}

\newtheorem{proposition}[theorem]{Proposition}
\newtheorem{definition}[theorem]{Definition}
\newtheorem{example}[theorem]{Example}
\newtheorem{remark}[theorem]{Remark}
\newtheorem{corollary}[theorem]{Corollary}

\def\QED{\quad\blackslug\lower 8.5pt\null}

\newcommand{\crazy}[2]{\displaystyle{\mathop{#1}_{#2}}
\vphantom{\displaystyle{#1}}}

\begin{document}

\begin{center}
{\Large \bf Maximum rank webs are not necessarily

\vspace*{1mm}

almost Grassmannizable}

\vspace*{5mm}

Vladislav V. Goldberg

\end{center}


\vspace*{3mm}

{\footnotesize \noindent {\em $2000$ Mathematics Subject
Classification}.  Primary 53A60
\newline

\textbf{Abstract} \textit{We present an example of a $6$-web $W
(6, 3, 2)$ of codimension two and of maximum rank on a
six-dimensional manifold which is not almost Grassmannizable.} }


\setcounter{section}{-1}

\setcounter{equation}{0}

\section{Introduction}
In 1984, during the Problem Session in the meeting on web geometry
at the Mathematisches Forschungsinstitut Oberwolfach,  Goldberg
posed the following problem (see \cite{S}):

\vspace*{2mm}

\textbf{6.} {\em Every $d$-web $W (d, n, r)$ of maximum $r$-rank
is almost Grassmannizable. Is it true or wrong?''}

Little in his paper \cite{L} related the Chern and Griffiths
approach in studying the rank problem and  Grassmannization and
algebraization problems for webs based on the presence of abelian
equations on the web with the Akivis and Goldberg approach based
on the notion of the almost Grassmannizable web. He wrote in the
introduction: ''The major purpose of this paper is to relate these
two approaches by showing that, in rough terms, if $d$ is large
enough (relative to $r, n$) so that $\pi (d, n, r) > 1$, then
every maximum rank web $W (d, n, r)$ is almost-Grassmannizable.
This answers (part of) a question posed by Goldberg at the 1984
Oberwolfach Conference on web geometry.''

Little proved in \cite{L} that
\textit{if $r > d (n-1) + 2$, then every maximum $r$-rank web
$W (d, n, r)$  is almost Grassmannizable.}

In \cite{L} Little also considered maximum $r$-rank webs $W (d, n,
r)$ with $d =  r (n-1) + 2$. In particular,  he proved that ``the
maximum 2-rank webs $W (6, 3, 2)$ are also almost
Grassmannizable'' (see his Example 1). It appears that this last
Little's result was wrong.

When Goldberg presented for publication the first version of his
paper \cite{G2}, where he considered the maximum rank webs $W (6,
3, 2)$, he  applied the above Little's result and deduced that
these kind of webs are always almost Grassmannizable.

A referee suspected that Little's result is incorrect and
provided a counterexample.
Goldberg recognized that the counterexample is correct
and communicated this counterexample to Little.
Little double-checked his proof and found
that the proof was incorrect. He discovered that the error in his proof
is in Corollary 3.6, and this happened because of rather tricky
a general position question. As a result, his Example 1, in which
he claimed that the webs $W (6, 3, 2)$ of maximum rank are almost Grassmannizable,
was incorrect.  The referee's counterexample was a counterexample to
Little's Corollary 3.6.

However, Little discussed the case
of webs $W (6, 3, 2)$ as a special case separately from the main
results which requires that the number of abelian equations be two
or more. The main result of his paper and his Example 2
giving an almost Grassmannizable web $W (r+2, 2, r)$ of maximum
rank one were correct.

In the revised version of the paper \cite{G2}, Goldberg
changed the reference to Little's paper \cite{L} by
the additional assumption that {\em the maximum rank webs
$W (d, n, r)$ in question are almost Grassmannizable.}

While looking recently at the paper \cite{G2}, the author of the
present paper recognized the importance of the referee
counterexample. This counterexample looks easy but it is not
trivial and was not known so far.

The purpose of this article is to present and study in detail the
above mentioned referee counterexample of a maximum 2-rank web $W
(6, 3, 2)$ which is  not almost Grassmannizable.

\section{Codimension Two Webs on a 6-Dimensional
Differentiable Manifold}

\textbf{1.} In an open domain $D$ of a differentiable manifold
$X^{6}$ of dimension six a 6-{\em web $W(6, 3, 2)$ of codimension
two}  is given by six codimension two foliations $X_\xi,
\linebreak \xi = 1, 2, 3, 4, 5, 6$, if the tangent 4-planes to the
leaves $V_\xi \subset X_\xi$ through a point in $D$ are in general
position.

Two webs $W(6,3,2)$ and $\widetilde{W}(6,3,2)$ are {\em
equivalent} to each other if there exists a local diffeomorphism
$\phi: D \rightarrow \widetilde{D}$ of their domains transferring
the foliations of $W$ into the foliations of $\widetilde{W}$.

Let $X_{\xi}, \;\; \xi = 1,\ldots, 6$, be six foliations of
parallel $4$-planes in an affine space $\mathbb{A}^{6}$ of
dimension six. Suppose that the $4$-planes of different foliations
are in general position. Such a $6$-web is called
\textit{parallel}. A web $W (6, 3, 2)$ which is equivalent to a
parallel web $W (6, 3, 2)$ is called \textit{parallelizable}.

The foliations $X_{\xi} ,\; \xi = 1, \ldots,6$, of the web
$W(6,3,2)$ can be given by six completely integrable systems of
Pfaffian equations
\begin{equation}\label{eq:1}
\crazy{\omega}{\xi}^i = 0, \;\;\;\; \xi = 1, \ldots ,6;
\;\;\;  i = 1, 2,
\end{equation}
where the forms
$\crazy{\omega}{\alpha}^i, \, \alpha = 1, 2, 3,$
are the basis forms of the manifold $X^{6}$ and
\begin{equation}\label{eq:2}
-\crazy{\omega}{4}^i =  \crazy{\omega}{1}^i
+ \crazy{\omega}{2}^i +  \crazy{\omega}{3}^i,
\end{equation}
\begin{equation}\label{eq:3}
- \crazy{\omega}{a}^i =
\crazy{\lambda}{a1}_j^i \;\crazy{\omega}{1}^j
+ \crazy{\lambda}{a2}_j^i \;\crazy{\omega}{2}^j
+ \crazy{\omega}{3}^i ,\;\;\;
a = 5, 6
\end{equation}
(see \cite{G1} for webs $W(d,2,r)$),
where the quantities
 $\crazy{\lambda}{a\hat{\alpha}}_j^i ,\; i,j = 1, 2$,
form an $(1,1)$-tensor for any $a = 5,6$ and $\hat{\alpha} = 1,
2$, and these four tensors $\crazy{\lambda}{a\hat{\alpha}}_j^i$
are distinct and satisfy some additional conditions. They are
called the {\it basis affinors} of a web $W(6,3,2)$ (cf.
\cite{G1}).

\vspace*{2mm} \textbf{2.}  We define now the notion of the almost
Grassmann structure. To this end, we need to define first the
notion of Segre cones (see \cite{AG}, Section 4). Consider the
Pl\"{u}cker mapping of the Grassmannian $G(2, 4)$ of $2$-planes of
a projective space $\mathbb{P}^{4}$ onto an algebraic manifold
$\Omega(2, 4)$ of dimension six of a projective space
$\mathbb{P}^9$. This mapping can be constructed by means of the
Grassmann coordinates of a $2$-plane $L$ in $P^{4}$  which are the
determinants of order three of the matrix
$$
\renewcommand{\arraystretch}{1.5}
\left(
\begin{matrix}
 x_1^1 & x_1^2 & x_1^3& x_1^{4} & x_1^{5} \cr
 x_2^1 & x_2^2 & x_2^3& x_2^{4} & x_2^{5} \cr
 x_3^1 & x_3^2 & x_3^3& x_3^{4} & x_3^{5} \cr
 \end{matrix} \right)
 \renewcommand{\arraystretch}{1}
$$
composed of the coordinates of the basis points $x_1, x_2, x_3$ of
the $2$-plane $L$. The Grassmann coordinates are connected by a
set of quadratic relations that define the  manifold $\Omega(2,
4)$ in the  space $P^9$ (see \cite{HP}, Chap. 7, \S 6). We will
say that this manifold carries the \textit{Grassmann structure}
and denote this manifold shortly by $\Omega$.

Let $L_1$ and $L_2$ be two $2$-planes in $\mathbb{P}^4$ meeting in
the straight line  $K$. They generate a linear pencil $S$ of
$2$-planes $\lambda L_1 + \mu L_2$. A rectilinear generator of the
manifold $\Omega$ corresponds to this pencil. All the $2$-planes
of the pencil $S$ belong to a $3$-plane $M$. This pencil, and
consequently the corresponding straight line in $\Omega$, is
completely determined by a pair $K$ and $M$, $K \subset M$.

Consider a bundle of $2$-planes, i.e., a set of all $2$-planes
passing through a fixed straight line $K$. On the manifold
$\Omega$, to this bundle there corresponds a two-dimensional plane
generator $\xi^2$.  On the other hand, on $\Omega$, to a family of
$2$-planes belonging to a fixed $3$-plane $M$, there corresponds a
three-dimensional plane generator $\eta^3$. Thus, the manifold
$\Omega$ carries two families of plane generators of dimensions
two and three, respectively.

If the straight line $K$ and the plane $M$ are incident,  $K
\subset M$, the plane generators $\xi^2$ and $\eta^3$ defined by
these planes, meet along a straight line. If they are not
incident, then the  generators $\xi^2$ and $\eta^3$  have no
common points.

Let us consider  a fixed $2$-plane $L$ in $\mathbb{P}^4$. It
contains a two-parameter family of straight lines $K$. Therefore,
the two-parameter family of generators $\xi^{2}$ passes through
the point $p \in  \Omega$ corresponding to $L$. On the other hand,
a one-parameter family of $3$-planes $M$ passes through the same
plane $L$. Consequently, a one-parameter family of generators
$\eta^{3}$ passes through the point $p$. Furthermore, any two
generators $\xi^{2}$ and $\eta^{3}$ passing through $p$ meet along
a straight line. It follows that all the plane generators
$\xi^{2}$ and $\eta^{3}$ passing through the point $p$, form a
cone whose projectivization is the \textit{Segre manifold} $S (1,
2)$ in the projective space $\mathbb{P}^5$. The Segre manifold
carries two families of plane generators of dimensions one and two
and can be considered as the projective embedding of the Cartesian
product of two projective spaces $\mathbb{P}^{1}$ and
$\mathbb{P}^{2}$ into the space $\mathbb{P}^{5}$:
$$
S (1, 2): \mathbb{P}^{1} \times \mathbb{P}^{2} \rightarrow
\mathbb{P}^{5}.
$$
The above described cone, whose projectivization is the Segre
manifold $S (1, 2)$, is called the \textit{Segre cone} and is
denoted by $C_{p} (2, 3)$. This cone is the intersection of the
manifold $\Omega$ and its tangent space $T_{p}(\Omega)$, whose
dimension is the same as that of $\Omega$, namely  six. In the
space $\mathbb{P}^{5}$, the set of all two-dimensional planes
intersecting a fixed $2$-plane $L$ along straight lines
corresponds to the cone $C_{p}(2, 3)$.

Therefore, with each point $p$ of the algebraic manifold $\Omega
\subset P^9$, there is connected the Segre cone $C_{p}(2, 3)$ with
vertex $p$ located in the tangent space $T_p(\Omega)$, and the
generators of this cone are generators of the manifold $\Omega$.

An \textit{almost Grassmann structure} on a manifold $X^6$ is a
smooth field of Segre cones (or their projectivizations) in the
tangent spaces $T X^6$.

Consider 4-subwebs [4, 1, 2, 3], [5, 1, 2, 3], and [6, 1, 2, 3]
defined by the first three foliations $X_1, X_2, X_3$ of a web $W
(6, 3, 2)$ and the foliations $X_4, X_5$, and $X_6$, respectively.
Each of these 4-subwebs defines an almost Grassmann structure. If
all  three  almost Grassmann structure so defined coincide, a web
$W (6, 3, 2)$ is said to be an \textit{almost Grassmann web}. Webs
$W (6, 3, 2)$ that are equivalent to an almost Grassmann web are
called \textit{almost Grassmannizable}. We denote such  webs by
$AGW(6, 3, 2)$. At each point all two-fold intersections of the
tangent spaces to the web leaves are 2-dimensional generators of
the Segre cone at that point.

For almost Grassmannizable webs $AGW(6,3,2)$, we have
\begin{equation}\label{eq:4}
\crazy{\lambda}{a\hat{\alpha}}_j^i =
\crazy{\lambda}{a\hat{\alpha}}\; \delta_j^i, \;\;\;a = 5,6;
\;\;\;\hat{\alpha} = 1,2
\end{equation}
(cf. \cite{G1}). Using equations (4), we can write  equations (3)
in the form
\begin{equation}\label{eq:5}
- \crazy{\omega}{a}^i =
\crazy{\lambda}{a1}\; \crazy{\omega}{1}^i
+ \crazy{\lambda}{a2}\; \crazy{\omega}{2}^i
+ \crazy{\omega}{3}^i,\;\;\;
a = 5, 6.
\end{equation}

The coefficients $\crazy{\lambda}{a1}$ and $\crazy{\lambda}{a2}$
in (4) satisfy certain inequalities which are implied by the fact
that the system of equations (1), (2), and (5) must be solvable
with respect to any six forms $\crazy{\omega}{\xi}^i, \,
\crazy{\omega}{\eta}^j, \, \crazy{\omega}{\zeta}^k, \; \xi, \eta,
\zeta = 1, \ldots, 6;\;\; \xi \neq \eta,  \zeta; \, \eta \neq
\zeta$ (see \cite{G2}).

The forms $\crazy{\omega}{\alpha}^i, \; \alpha = 1, 2, 3; \, i =
1, 2,$ satisfy the following structure equations:
\begin{equation}\label{eq:6}
d\crazy{\omega}{\alpha}^{i} - \crazy{\omega}{\alpha}^{j} \wedge
\omega_{j}^{i} =  \sum_{\beta\neq\alpha}
\crazy{a}{\alpha\beta}^{i}_{jk} \; \crazy{\omega}{\alpha}^{j}
\wedge \crazy{\omega}{\beta}^{k}
\end{equation}
 (see \cite{G1}), where the quantities $\crazy{a}{\alpha\beta}^{i}_{jk}$
form the \textit{torsion tensor} of the web $AGW(6,3,2)$.

On the manifold $X^6$, a web $W(6,3,2)$ defines an affine connection
$\gamma$ which is determined by the forms
$\omega^\xi = \{\crazy{\omega}{\alpha}^i\},
 \;\xi, \eta = 1,2,3,4,5,6;\; \alpha = 1,2,3;\; i = 1, 2$, and
$$
\renewcommand{\arraystretch}{1.5}
  \left(
  \begin{array}{lll}
        (\omega^i_j) & O_{2\times 2} & O_{2\times 2} \\
        O_{2\times 2}  & (\omega^i_j)  & O_{2\times 2} \\
        O_{2\times 2}  & O_{2\times 2}  &  (\omega^i_j)
\end{array}
\right)
\renewcommand{\arraystretch}{1}
$$
(see \cite{G1}), where $(\omega_i^j) =
  \left(
  \begin{array}{ll}
        \omega^1_1 & \omega^2_1  \\
        \omega^1_2  & \omega^2_2 \\
\end{array}
\right)
$ and $O_{2\times 2} =
 \left(
  \begin{array}{ll}
        0 & 0  \\
        0 & 0 \\
\end{array}
\right).
$

\vspace*{2mm}

\textbf{3.} Suppose that the leaves of the  $\xi$th foliation of a
web $W(6,3,2)$ are given as level sets of functions $u_{\xi}^i(x)$
:
$$
 u_{\xi}^i (x) = \mbox{const.},\;\;\;\;\; \xi = 1, \ldots ,6.
$$
The functions $u_{\xi}^i (x)$ are defined up to a local
diffeomorphism in the space of $u_{\xi}^i(x)$.

An exterior 2-equation of the form
$$
 \sum_{\xi = 1}^6 f_{\xi} (u_{\xi}^j)\; du_{\xi}^1 \wedge
 du_{\xi}^2 = 0, \;\;\;\; j =  1, 2,
$$
is said to be an \textit{abelian $2$-equation}. The maximum number
$R_2$ of linearly independent abelian 2-equations admitted by a
web $W(6,3,2)$ is called the \textit{$2$-rank} of the web
$W(6,3,2)$.

It follows from the definition that the coefficients
$f_{\xi}$ are constant on the leaves of the  $\xi$th
foliation of $W(6,3,2)$.

If there exists an upper bound $\pi_2 (6,3,2)$ of $R_2$,
then $R_2 \leq \pi_2(6,3,2)$.
Chern and Griffiths \cite{CG1} found the general formula
for $\pi (d, n, r)$. It follows from their formula
that $\pi_2(6,3,2) = 1$.

It was proved in \cite{G2} that the condition
\begin{equation}\label{eq:7}
 \crazy{\lambda}{51} \crazy{\lambda}{62} \,(1 - \crazy{\lambda}{52} -  \crazy{\lambda}{61})
 = \crazy{\lambda}{52} \crazy{\lambda}{61} \,(1 - \crazy{\lambda}{51} - \crazy{\lambda}{62})
 \end{equation}
 is a necessary condition for a web $AGW (6, 3, 2)$ to be of maximum
 rank one.

\section{Referee's Counterexample}

\textbf{4.} Suppose that $(x^1, \, x^2, \, x^3, \, y_4, \, y_5, \,
y_6)$ are coordinates in $\mathbb{R}^6$ and $A$ is any nonsingular
$3 \times 3$ matrix of real numbers. Define the linear functions
$x^4, \, x^5, \, x^6, \, y_1, \, y_2, \, y_3$  on $\mathbb{R}^6$
as follows:
$$
(x^4, x^5, x^6) = (x^1, x^2, x^3)\, A \;\; \mbox{and} \;\; ( y_1,
y_2, y_3) = - (y_4, y_5, y_6)\, A^T.
$$
If
$$
\renewcommand{\arraystretch}{1.5}
A = \left(
\begin{array}{lll}
a_1^1 & a_1^2 & a_1^3 \\
             a_2^1 & a_2^2 & a_2^3 \\
             a_3^1 & a_3^2 & a_3^3
\end{array}
\right),
\renewcommand{\arraystretch}{1}
$$
then
\begin{equation}\label{eq:8}
x^{3+\alpha} = a_\beta^\alpha \,x^\beta, \;\; \alpha, \beta = 1,
2, 3,
 \end{equation}
and
\begin{equation}\label{eq:9}
y_\alpha =  - a^\beta_\alpha \,y_{3 + \beta}, \;\; \alpha, \beta =
1, 2, 3.
 \end{equation}
Equation (9) can be written as
\begin{equation}\label{eq:10}
y_{3+\alpha} = -  b_\alpha^\beta \,y_\beta,
 \end{equation}
 where $B = (b_\alpha^\beta)$ is the inverse matrix of the matrix
 $A$.

Differentiation of (8) and (9) gives
\begin{equation}\label{eq:11}
dx^{3+\alpha} = a_\beta^\alpha \,dx^\beta,
 \end{equation}
and
\begin{equation}\label{eq:12}
dy_\alpha =  - a^\beta_\alpha \, dy_{3 + \beta}.
 \end{equation}
Using (11) and (12), we can easily find that
\begin{equation}\label{eq:13}
d x^\xi \wedge dy_\xi = 0, \;\; \xi = 1, 2, 3, 4, 5, 6.
 \end{equation}

Hence equations
$$
dx^\xi = 0, \;\; dy_\xi = 0, \;\; \xi =  1, 2, 3, 4, 5, 6,
$$
define a set of foliations $X_\xi$ of codimension two on
$\mathbb{R}^6$, i.e., a 6-web $W (6, 3, 2)$ admitting an abelian
2-equation (13).

First, we prove  that the 6-web we have constructed is
parallelizable.

\begin{proposition}
The $6$-web defined by equations $(8)$ and $(10)$ is
parallelizable.
\end{proposition}

\begin{proof}
In fact, each of the foliations $X_\xi$ of the constructed 6-web
are defined by equations
$$
dx^\xi = 0, \;\; dy_\xi = 0, \;\; \xi =  1, 2, 3, 4, 5, 6,
$$
where the index $\xi$ is fixed, or
$$
x^\xi = \text{const.}, \;\; y_\xi = \text{const.}
$$
Each of these two equations represents a foliation of parallel
hyperplanes in $\mathbb{R}^6$, and the leaves of $X_\xi$ are
parallel 4-planes that are intersections of parallel hyperplanes
of the above two foliations of hyperplanes.
\end{proof}

\textbf{5.} According to \cite{CG1}, a 6-web can admit at most one
abelian 2-equation. For a generic matrix $A$, the web normals $d
x^1 \wedge dy_1, d x^2 \wedge dy_2, d x^3 \wedge dy_3, d x^4
\wedge dy_4, d x^5 \wedge dy_5$, and $d x^6 \wedge dy_6$ are, in
some sense, in the most general position which allows there to be
an abelian 2-equation (13). So, we constructed a 6-web $W (6, 3,
2)$ of maximum 2-rank. Note that because the matrix $A$ is
generic, its entries (which are constants) are not connected by
any relation, i.e., there is a 9-parameter family of such
matrices.

We  prove now that the 6-web constructed above is not almost
Grassmannizable.

\begin{theorem}
The $6$-web defined by equations $(8)$ and $(10)$ is not almost
Grassmannizable.
\end{theorem}

\begin{proof}
To prove the theorem, first we  reduce equations (11) and (12) to
the form (2) and (3).

Denote
\begin{equation}\label{eq:14}
\crazy{\omega}{\alpha}^1 = dx^\alpha, \;\;
\crazy{\omega}{\alpha}^2 =dy_\alpha = - a_\alpha^\beta\,
dy_{3+\beta}.
 \end{equation}

Solving equations (12), we find that
\begin{equation}\label{eq:15}
-dy_{3+\alpha} =  \sum_\beta b_\alpha^\beta\, dy_\beta,
 \end{equation}
 where $B = (b_\alpha^\beta)$ is the inverse matrix of the matrix
 $A$.

Denote
\begin{equation}\label{eq:16}
 \crazy{\omega}{4}^1 =  dx^4, \;\;
 \crazy{\omega}{4}^2 = dy_4.
 \end{equation}
Define
\begin{equation}\label{eq:17}
- \crazy{\overline{\omega}}{\alpha}^1 = a_\alpha^1 \,
\crazy{\omega}{\alpha}^1,\;\;
 \crazy{\overline{\omega}}{\alpha}^2 = b_1^\alpha \, \crazy{\omega}{\alpha}^2
\;\; (\mbox{no summation over $\alpha$}).
 \end{equation}
If we suppress the bar over $\crazy{\omega}{\alpha}^i$, then
equations (16) take the form (2):
\begin{equation}\label{eq:18}
- \crazy{\omega}{4}^i =  \crazy{\omega}{1}^i + \crazy{\omega}{2}^i
+ \crazy{\omega}{3}^i.
 \end{equation}

It follows from (11), (15), and (17) that
\begin{equation}\label{eq:19}
\renewcommand{\arraystretch}{2.7}
\left\{
\begin{array}{ll}
- dx^5 = \displaystyle  \frac{a_1^2}{a^1_1}\; \crazy{\omega}{1}^1
+ \frac{a_2^2}{a^1_2}\; \crazy{\omega}{2}^1 + \frac{a_3^2}{a^1_3}
\;\crazy{\omega}{3}^1, & - dx^6 = \displaystyle
\frac{a_1^3}{a^1_1} \; \crazy{\omega}{1}^1 + \frac{a_2^3}{a^1_2}
\;\crazy{\omega}{2}^1
+ \frac{a_3^3}{a^1_3} \;\crazy{\omega}{3}^1, \\
- dy_5 = \displaystyle \frac{b_2^1}{b^1_1} \;\crazy{\omega}{1}^2 +
\frac{b_2^2}{b^2_1} \; \crazy{\omega}{2}^2 + \frac{b_2^3}{b^3_1}
\;\crazy{\omega}{3}^2, & - dy_6 = \displaystyle
\frac{b_3^1}{b^1_1} \;\crazy{\omega}{1}^2 + \frac{b_3^2}{b^2_1}\;
\crazy{\omega}{2}^2 + \frac{b_3^3}{b^3_1} \; \crazy{\omega}{3}^2.
\end{array}
\right.
\renewcommand{\arraystretch}{1}
 \end{equation}
Equations (19) can be written as
\begin{equation}\label{eq:20}
\renewcommand{\arraystretch}{2.7}
\left\{
\begin{array}{ll}
- \displaystyle  \frac{a_3^1}{a^2_3} \; dx^5 = \displaystyle
\frac{a_1^2 a_3^1}{a^1_1 a_3^2} \; \crazy{\omega}{1}^1 +
\frac{a_2^2 a_3^1}{a^1_2 a_3^2} \; \crazy{\omega}{2}^1 + \crazy{\omega}{3}^1, \\
- \displaystyle  \frac{a_3^1}{a^3_3} \; dx^6 = \displaystyle
\frac{a_1^3 a_3^1}{a^1_1 a_3^3} \; \crazy{\omega}{1}^1
+ \frac{a_2^3 a_3^1}{a^1_2 a_3^3} \; \crazy{\omega}{2}^1 +  \crazy{\omega}{3}^1, \\
- \displaystyle \frac{b_1^3}{b^3_2} \;dy_5 = \displaystyle
\frac{b_2^1 b^3_1}{b^1_1 b^3_2}\; \crazy{\omega}{1}^2 +
\frac{b_2^2 b^3_1}{b^2_1 b^3_2}\; \crazy{\omega}{2}^2 +
\crazy{\omega}{3}^2,
\\ - \displaystyle  \frac{b_1^3}{b^3_3} \;dy_6 = \displaystyle
\frac{b_3^1 b^3_1}{b^1_1 b^3_3} \; \crazy{\omega}{1}^2 +
\frac{b_3^2 b^3_1}{b^2_1 b^3_3} \; \crazy{\omega}{2}^2 +
\crazy{\omega}{3}^2.
\end{array}
\right.
\renewcommand{\arraystretch}{1}
 \end{equation}

Define
\begin{equation}\label{eq:21}
 \crazy{\omega}{5}^1 = \displaystyle \frac{a_3^1}{a^2_3} \;
 dx^5, \;\; \crazy{\omega}{5}^2 = \displaystyle  \frac{b_1^3}{b^3_2}
 \; dy_5, \;\;
 \crazy{\omega}{6}^1 = \displaystyle  \frac{a_3^1}{a^3_3}
 \; dx^6, \;\;  \crazy{\omega}{6}^2 = \displaystyle  \frac{b_1^3}{b^3_3}
 \; dy_6.
 \end{equation}
Then equations (20) take the form
\begin{equation}\label{eq:22}
\renewcommand{\arraystretch}{2.7}
\left\{
\begin{array}{ll}
- \crazy{\omega}{5}^1 = \displaystyle \frac{a_1^2 a_3^1}{a^1_1
a_3^2} \; \crazy{\omega}{1}^1 + \frac{a_2^2 a_3^1}{a^1_2 a_3^2} \;
\crazy{\omega}{2}^1 + \crazy{\omega}{3}^1, & - \crazy{\omega}{6}^1
= \displaystyle \frac{a_1^3 a_3^1}{a^1_1 a_3^3}\;
\crazy{\omega}{1}^1
+ \frac{a_2^3 a_3^1}{a^1_2 a_3^3}\; \crazy{\omega}{2}^1 +  \crazy{\omega}{3}^1, \\
- \crazy{\omega}{5}^2 = \displaystyle \frac{b_2^1 b^3_1}{b^1_1
b^3_2} \; \crazy{\omega}{1}^2 + \frac{b_2^2 b^3_1}{b^2_1 b^3_2} \;
\crazy{\omega}{2}^2 +  \crazy{\omega}{3}^2, & -
\crazy{\omega}{6}^2 = \displaystyle \frac{b_3^1 b^3_1}{b^1_1
b^3_3} \; \crazy{\omega}{1}^2 + \frac{b_3^2 b^3_1}{b^2_1 b^3_3} \;
\crazy{\omega}{2}^2 + \crazy{\omega}{3}^2.
\end{array}
\right.
\renewcommand{\arraystretch}{1}
 \end{equation}

Comparing equations (22) with equations (3), we find that
\begin{equation}\label{eq:23}
\renewcommand{\arraystretch}{2.7}
\left\{
\begin{array}{llll}
\crazy{\lambda}{51}_i^1 = \frac{a_1^2 a_3^1}{a^1_1 a_3^2} \;
\delta_i^1, & \crazy{\lambda}{52}_i^1 = \frac{a_2^2 a_3^1}{a^2_1
a_3^2} \; \delta_i^1,& \crazy{\lambda}{51}_i^2 = \frac{b_2^1
b_1^3}{b^1_1 b_2^3} \; \delta_i^2,& \crazy{\lambda}{52}_i^2 =
\frac{b_2^2 b_1^3}{b^1_1 b_2^3} \; \delta_i^2,
\\
\crazy{\lambda}{61}_i^1 = \frac{a_1^3 a_3^1}{a^1_1 a_3^3} \;
\delta_i^1, & \crazy{\lambda}{62}_i^1 = \frac{a_2^3 a_3^1}{a^2_1
a_3^3} \; \delta_i^1,& \crazy{\lambda}{61}_i^2 = \frac{b_3^1
b_1^3}{b^1_1 b_3^3} \; \delta_i^2,& \crazy{\lambda}{62}_i^2 =
\frac{b_3^2 b_1^3}{b^1_2 b_3^3} \; \delta_i^2.
\end{array}
\right.
\renewcommand{\arraystretch}{1}
\end{equation}

These equations show that \textit{a web $W (6, 3, 2)$
of maximum $2$-rank defined by equations $(8)$ and $(10)$
has constant basis affinors.}

Next suppose that our 6-web  under consideration is almost
Grassmannizable. This will be the case if and only if equations
(22) have the form (5) or the functions
$\crazy{\lambda}{a\hat{\alpha}}_j^i$ have the form (4). Comparing
(4) and (23), we find that our web is almost Grassmannizable if
and only if the following four conditions hold:
\begin{equation}\label{eq:24}
\renewcommand{\arraystretch}{2.7}
\left\{
\begin{array}{ll}
\crazy{\lambda}{51} = \frac{a_1^2 a_3^1}{a^1_1 a_3^2} = \frac{b_2^1 b_1^3}{b^1_1 b_2^3},
&
\crazy{\lambda}{52} = \frac{a_2^2 a_3^1}{a^2_1 a_3^2} = \frac{b_2^2 b_1^3}{b^1_1 b_2^3},
\\
\crazy{\lambda}{61} = \frac{a_1^3 a_3^1}{a^1_1 a_3^3} = \frac{b_3^1 b_1^3}{b^1_1 b_3^3},
&
\crazy{\lambda}{62} = \frac{a_2^3 a_3^1}{a^2_1 a_3^3} =  \frac{b_3^2 b_1^3}{b^1_2 b_3^3}.
\end{array}
\right.
\renewcommand{\arraystretch}{1}
\end{equation}

A straightforward calculation shows that these four conditions are
equivalent to the condition which can be written in three equivalent
forms:
\begin{equation}\label{eq:25}
\renewcommand{\arraystretch}{1.5}
\left|
\begin{array}{ccc}
a_1^1 & a_1^2 & a_1^1 a_1^2 a_2^3 a_3^3\\
a_2^1 & a_2^2 & a_2^1 a_2^2 a_1^3 a_3^3\\
a_3^1 & a_3^2 & a_1^3 a_2^3 a_3^1 a_3^2
\end{array}
\right| = 0,
\renewcommand{\arraystretch}{1}
\end{equation}
or
\begin{equation}\label{eq:26}
\renewcommand{\arraystretch}{1.5}
\left|
\begin{array}{ccc}
a_1^1 & a_1^1 a_1^3 a_2^2 a_3^2 & a_1^3\\
a_2^1 & a_1^2 a_2^1 a_2^3 a_3^2 & a_2^3\\
a_3^1 & a_1^2 a_2^2 a_2^2 a_3^1 & a_3^3
\end{array}
\right| = 0,
\renewcommand{\arraystretch}{1}
\end{equation}
or
\begin{equation}\label{eq:27}
\renewcommand{\arraystretch}{1.5}
\left|
\begin{array}{ccc}
a_1^1 a_1^3 a_2^2 a_3^2 & a_1^2&  a_1^3\\
a_1^2 a_2^1 a_2^3 a_3^2 & a_2^2 & a_2^3\\
a_1^2 a_2^2 a_2^2 a_3^1 & a_3^2 & a_3^3
\end{array}
\right| = 0.
\renewcommand{\arraystretch}{1}
\end{equation}

We came to the contradiction, since equation (25) (or (26), or
(27)) shows that the matrix $A$ is not the most generic. In fact,
the matrices $A$ satisfying  equation (25) (or (26), or (27)) form
an 8-parameter subfamily in the 9-parameter family of the most
generic matrices $A$.
\end{proof}

\emph{Remark.} As we noted earlier, if an almost Grassmannizable
web $AGW (6, 3, 2)$ is of maximum $2$-rank (i.e., it admits one
abelian 2-equation), then it is necessary that condition (7)
holds. Since our web admits an abelian 2-equation (13), condition
(7) should be satisfied identically. A simple calculation shows
that in fact this is true: by (24), condition (7) is reduced to
the same equation (25) (or (26), or (27)) which holds for an
almost Grassmannizable web $AGW (6, 3, 2)$.

\vspace*{2mm}

\textbf{6.} We present now another (analytic) proof of Proposition
1.

First, we remind that a web $W(d,n,r)$ is said to be
\textit{parallelizable} if it is equivalent to a web $W(d,n,r)$
formed by $d$  foliations of parallel $(n-1)r$-dimensional planes
of $(nr)$-dimensional affine space $A^{nr}$.

The following general criteria of parallelizability is valid:
\emph{ A web $W(d,n,r)$ is parallelizable if and only all its
basis affinors $\crazy{\lambda}{a \sigma}^i_j,\, a = n+2, \ldots ,
d;\, \sigma = 1, \ldots, n - 1$, are covariantly constant on
$X^{nr}$ in an affine connection $\gamma$ and its $(n+1)$-subweb
$[1,2, \ldots, n]$ defined by the foliations $X_\alpha, \, \alpha
= 1, 2, \ldots , n$,  is parallelizable. }

The proof of this theorem is similar to that of Theorem
7.2.2 in \cite{G1} (p. 311) proved there for webs $W (4, 2, r)$.

\begin{proof}
In fact, taking exterior derivatives of equations (14), we find that
\begin{equation}\label{eq:28}
d\crazy{\omega}{\alpha}^i = 0, \;\;\;\; \alpha = 1, 2, 3; \; i =
1, 2.
\end{equation}
Comparing equations (28) with the structure equations (6) of a
general web $W (6, 3, 2)$, we see that
 $$
\omega_{j}^{i} = 0, \;\; \crazy{a}{\alpha\beta}^{i}_{jk} = 0.
 $$

Since the basis affinors (23) of a web  defined by equations (8) and (10)
are constant and $\omega_{j}^{i} = 0$, it follows that
$$
\nabla \crazy{\lambda}{a \sigma}^i_j = d \crazy{\lambda}{a
\sigma}^i_j - \crazy{\lambda}{a \sigma}^i_k \,\omega^k_j +
\crazy{\lambda}{a \sigma}^k_j \,\omega_k^i = 0,
$$
i.e., the basis affinors $\crazy{\lambda}{a \sigma}^i_j$ are
covariantly constant on $\mathbb{R}^6$. By Theorem 1, this implies
that the webs $W (6, 3, 2)$ of maximum $2$-rank defined by
equations (8) and (10) are  parallelizable.
 \hfill \end{proof}

The following corollary immediately follows from Proposition 1 and
Theorem 2:

\begin{corollary}
Equations   $(8)$ and $(10)$ define a $9$-parameter family of
parallelizable, not almost Grassmannizable webs of maximum
$2$-rank.
\end{corollary}

\vspace*{2mm}

\textbf{7.} As we indicated earlier, there is a 9-parameter family
of  parallelizable, not almost Grassmannizable webs  of maximum
$2$-rank. Such webs can be obtained by choosing a nonsingular
matrix $A$ with the only condition that the determinant on the
left-hand side of equation (25) (or (26) or (27)) does not vanish.

For any $a = 1, 2, 3, 4, 5, 6, 7, 8$,
it is easy to find $a$-parameter subfamilies
of this  9-parameter family of webs.

For example, an 8-parameter family $B_8$ of such webs is defined
by the matrices $A$ with $a_1^3 = 0$; a 7-parameter family $B_7$
of such webs is defined by the matrices $A$ with $a_1^2 = a_1^3 =
0$; and a 6-parameter family $B_6$ of such webs is defined by the
matrices $A$ with $a_1^3 = a_2^1 = a_3^2 = 0$. In each of these
cases, it is easy to check that the determinant on the left-hand
side of (27) does not vanish.

For each of the subfamilies  $B_8, B_7$, and  $B_6$, we can define
$a$-parameter subfamilies for any $a = 1, \ldots , 7,$ by relating
some of the remaining nonvanishing entries of the matrix $A$.

We conclude by indicating three concrete examples of webs
belonging respectively to the subfamilies  $B_8, B_7$, and  $B_6$.

\vspace*{2mm}

\textbf{Example 1.} Consider
$$
A =
\left(
\begin{array}{lll}
1 & 1 & 0 \\
             1 & 1 & 1 \\
             1 & 2 & 1
\end{array}
\right).
$$
It is easy to see that this matrix is nonsingular, and the
determinant on the left-hand side of (27) is equal to 2.

A simple computation by means of equations (8) and (10) gives the
following closed form equations of this web:
\begin{equation}\label{eq:29}
\renewcommand{\arraystretch}{1.5}
\left\{
\begin{array}{lllll}
x^4 = x^1 + &\!\!\!\! x^2 + x^3,  & y_4 = - &\!\!\!\! y_1 - &\!\!\!\! y_2 +
y_3,\\
x^5 = x^1 + &\!\!\!\! x^2 + 2x^3, & y_5 =   &\!\!\!\!       &\!\!\!\!  y_2 -y_3,\\
x^6 =       &\!\!\!\! x^2 + x^3,  & y_6 =  &\!\!\!\! y_1 -  &\!\!\!\! y_2.
\end{array}
\right.
\renewcommand{\arraystretch}{1}
\end{equation}

\textbf{Example 2.} Consider
$$
\renewcommand{\arraystretch}{1.3}
A = \left(
\begin{array}{lll}
1 & 1 & 0 \\
             0 & 1 & 1 \\
             1 & 1 & 1
\end{array}
\right).
\renewcommand{\arraystretch}{1}
$$
It is easy to see that this matrix is nonsingular, and the
determinant on the left-hand side of (27) is equal to 1.

A simple computation by means of equations (8) and (10) gives the
following closed form equations of this web:
\begin{equation}\label{eq:30}
\renewcommand{\arraystretch}{1.5}
\left\{
\begin{array}{lllllll}
x^4 = x^1   &\!\!\!\!         &\!\!\!\! +x^3,  & y_4 =  &\!\!\!\!       &\!\!\!\!y_2  &\!\!\!\!- y_3,\\
x^5 = x^1 + &\!\!\!\! x^2,     &\!\!\!\!      & y_5 =  -&\!\!\!\!y_1 -  &\!\!\!\!y_2 &\!\!\!\! + y_3,\\
x^6 =       &\!\!\!\! x^2     &\!\!\!\! +x^3, & y_6 =  -&\!\!\!\!y^1    &\!\!\!\!    &\!\!\!\!+ y_3.
\end{array}
\right.
\renewcommand{\arraystretch}{1}
\end{equation}

\textbf{Example 3.} Consider
$$
\renewcommand{\arraystretch}{1.3}
A = \left(
\begin{array}{lll}
1 & 1 & 0 \\
             0 & 1 & 1 \\
             1 & 0 & 1
\end{array}
\right).
\renewcommand{\arraystretch}{1}
$$
It is easy to see that this matrix is nonsingular, and the
determinant on the left-hand side of (27) is equal to $-1$.

A simple computation by means of equations (8) and (10) gives the
following closed form equations of this web:
\begin{equation}\label{eq:31}
\renewcommand{\arraystretch}{1.5}
\left\{
\begin{array}{lllll}
x^4 = x^1   &\!\!\!\!       +  &\!\!\!\!x^3,  & y_4 =  \frac{1}{2} (-&\!\!\!\!y_1 + y_2 + y_3),\\
x^5 = x^1 + &\!\!\!\! x^2,     &\!\!\!\!      & y_5 =  \frac{1}{2} (-&\!\!\!\!y_1 - y_2 + y_3),\\
x^6 =       &\!\!\!\! x^2 +    &\!\!\!\! x^3, & y_6 =  \frac{1}{2} ( &\!\!\!\!y^1 - y_2 - y_3).
\end{array}
\right.
\renewcommand{\arraystretch}{1}
\end{equation}

\vspace*{2mm}

\textbf{8.} Finally note  that if in equations (8) and (10) the
indices $\alpha$ and $\beta$ take values $1, 2, \ldots , n$, and
$A$ is an arbitrary $n \times n$ matrix, then these equations
produce a $(2n)$-web $W (2n, n, 2)$ which is {\em a
parallelizable, not almost Grassmannizable $(2n)$-web of maximum
$2$-rank one}.

\makeatletter \renewcommand{\@biblabel}[1]{\hfill#1.}\makeatother
\bibliographystyle{amsplain}

\vspace*{10mm}

\noindent {\em Author's address}:\\

\noindent
V.~V. Goldberg\\
Department of Mathematical Sciences\\
New Jersey Institute of Technology \\
 University Heights \\
  Newark, N.J. 07102, U.S.A. \\
  \\
 E-mail address:
 vlgold@m.njit.edu
\end{document}